\documentclass{article}

\usepackage{PRIMEarxiv}

\usepackage[utf8]{inputenc} 
\usepackage[T1]{fontenc}    
\usepackage{hyperref}       
\usepackage{url}            
\usepackage{booktabs}       
\usepackage{amsfonts}       
\usepackage{nicefrac}       
\usepackage{microtype}      
\usepackage{lipsum}
\usepackage{fancyhdr}       
\usepackage{graphicx}       
\graphicspath{{media/}}     
\usepackage[bold,full]{complexity}

\title{Hard Proofs and Good Reasons}

\author{
   Simon DeDeo \\
   Carnegie Mellon University, Pittsburgh, PA 15213, USA \\
    \& the Santa Fe Institute, Santa Fe, NM 87501, USA \\
   \texttt{sdedeo@andrew.cmu.edu}
}

\begin{document}
\maketitle

\begin{abstract}
Practicing mathematicians often assume that mathematical claims, when they are true, have good reasons to be true. Such a state of affairs is ``unreasonable'', in Wigner's sense, because basic results in computational complexity suggest that there are a large number of theorems that have only exponentially-long proofs, and such proofs can not serve as good reasons for the truths of what they establish. Either mathematicians are adept at encountering only the reasonable truths, or what mathematicians take to be good reasons do not always lead to equivalently good proofs. Both resolutions raise new problems: either, how it is that we come to care about the reasonable truths before we have any inkling of how they might be proved, or why there should be good reasons, beyond those of deductive proof, for the truth of mathematical statements. Taking this dilemma seriously provides a new way to make sense of the unstable ontologies found in contemporary mathematics, and new ways to understand how non-human, but intelligent, systems might found new mathematics on inhuman ``alien'' lemmas.
\end{abstract}

\vspace{0.5cm}

The ``unreasonable'' effectiveness of mathematics in the sciences~\cite{wigner} is a topic of perennial fascination~\cite{coly}: in a world of baffling complexity, we find that core laws take a mathematically simple form. In Steiner's formulation of the problem~\cite{steiner2009applicability}, the laws of the universe have an unexpectedly user-friendly structure, where the guidance provided by mathematical taste not only picks out useful theories, but has its own generative quality: the cognitive appeal of a piece of mathematics anticipates structures that find unexpected and generative use in yet-to-be encountered physics.

We seem to see, in other words, an anthropocentric reason-to-things that is at odds with an evolutionary, or Copernican, view where human reason arrives as a late-comer to the universe. There are many explanations for this puzzling ``nomological harmony''~\cite{nomological}, ranging from dismissal~\cite{dismissal} or denial~\cite{bob,ab_ex} to the metaphysically elaborate, platonistic accounts by physicists such as Penrose~\cite{penrose1999emperor} and Tegmark~\cite{tegmark2008mathematical}.

The goal of this paper is to draw attention to a similar, less remarked-upon, phenomenon that occurs within mathematics itself. Mathematicians, like those in the natural sciences, also interact with a generative, non-human phenomenon: the space of deductively-valid proofs. They do so, necessarily, using tools provided by cultural evolution such as explanation-making~\cite{chater2016under}, argument~\cite{mercier2011humans} and dialectic~\cite{novaes2020dialogical}. Whatever these tools evolved to do---convince others to share one's beliefs, say, or navigate the world of ``medium-sized dry goods''~\cite{austin1964sense}---they seem to have escaped those roots to enable us to discover the deductive proofs of truly strange claims well beyond the training data provided in cultural history.

The point here is not that mathematicians contemplate mysteries far beyond the ordinary, nor that they convince their colleagues they have done so. There have always been mystics, and there have always been shared delusions. What is strange is that mathematicians seem to do this in a way that enables them to produce deductive proofs of their discoveries. Mathematicians may have their doubts about the rigor of this or that claim, but (to take one example) the successful computer-verification of elaborate and unintuitive constructions such as ``liquid tensors'' ~\cite{scholze2022liquid} is taken by most mathematicians to validate a norm of their activity.

Feeling that something is unreasonable, of course, does not prove it so, and the purpose of this paper is to show why this state of affairs is, in the folk picture, quite unreasonable given key theorems about deduction and proof-making from computational complexity. To sharpen the point, I will draw upon an account of doing mathematics provided by Timothy Gowers~\cite{gowers2023makes} in his recent paper ``What makes mathematicians believe unproved mathematical statements''. The goal here is not to take Gowers' testimony as truth, but rather to draw attention to the strangeness of what doing so would entail, and the ways in which it is in tension with the metamathematics we have from the theoretical computer sciences.

Gowers's insights are particularly valuable for those interested in the philosophy and psychology of mathematics because he presents his results as empirical generalizations of the heuristics that he, and other mathematicians, use on a daily basis. Asking exactly how and why an (apparent) anthropocentrism arises in mathematics gives us new ways to look at the practice of mathematics and to make sense of the testimony of mathematicans working at the cutting edge of contemporary mathematics; new ways to think about how artificial intelligence may change mathematics; and counter-intuitive possibilities such as the existence of yet-to-be discovered, fertile fields of mathematics based on novel axioms that are essentially impossible to found in the systems we have to hand.

\section{Believing the to-be-proven}

A great deal of work has been devoted to the general question of how mathematicians find proofs. Indeed, investigations of the search problem are one of the origin points of the entire field of computational complexity. Gowers (and I) deal with a more focused question: how we come to believe mathematical claims, to some degree of confidence, before we have the proofs in hand. This is not a small question, however, because trying to predict whether or not a mathematical statement is true is, as Gowers out, more than just a parlour game of polling the room for beliefs about famous conjectures. Forming intuitions about the likelihood that this or that statement is true is a core part of going about trying to prove things, because in the process of constructing a proof, a mathematicians is consistently faced with the question of what claims to try to establish next. Those decisions are guided in part by what (proto-)lemmas one thinks are likely to be true: if one is trying to prove $B$, and notices that $A$ implies $B$, the next question may well be whether or not $A$ is likely to be true.

Gowers presents a number of case studies to illustrate two key principles for how mathematicians go about conjecturing the truth of a mathematical statement: the {\bf no-coincidence principle} (``if an apparently outrageous coincidence happens in mathematics, then there is a reason for it'') and the {\bf no-miracles principle} (``if an apparent miracle happens in mathematics, then there is a reason for it''). The distinction between the two is whether or not the fact in question is unexpected on broadly probabilistic grounds\footnote{For a recent review, and example of arguments in this style, see Ref.~\cite{christiano2022formalizing}, and Ref.~\cite{taylor} for a more elaborate formulation of the interpersonal case.} (in which case it is a ``coincidence''), or on more general expectations about how the world works (in which case it is a ``miracle'').

To see how these principles play out in practice, consider the example of the (strong) Goldbach Conjecture (GC), which says that every even number larger than two is the sum of two primes. As Gowers points out, the GC has been verified, by computer, for even integers up to $4\times10^{18}$: for all of them, so far, the computer can find the required pair of primes. In addition to this raw data, there is good evidence that the prime numbers---treated empirically, in the way a scientist might study a natural phenomenon---obey statistical laws that fit the (now proven) weak Goldbach Conjecture, and imply the strong conjecture as well.

From this, Gowers argues, we learn two things. First, if Goldbach's Conjecture is false, it must be an ``outrageous coincidence'': prime numbers would have to fortuitously arrange themselves in such a way as to produce an outcome that was statistically surprising as a matter of Bayesian evidence. Second, that ``outrageous coincidence'' is unlikely to have a ``reason'', because such a reason would likely have already intervened to upset the relationship between the weak-GC and the statistical evidence that accumulated in its favor, or would be otherwise easy to find. ``No coincidence without reason'' (the no-coincidence principle) together with ``no reason'' (implied by experience and intuition), implies no coincidence, and thus that we ought to believe in the truth of the strong Goldbach Conjecture.

\section{Proofs, Reasons, and the No-Coincidence Principle}

Gowers' argument hinges on the notion of something having a reason to be true, and it is worth examining more closely what he means by this. Any proof can be seen as providing a reason why its claim is true, simply as a matter of deduction. If I ask, for example, why I have trouble dividing my class of 83 students into smaller groups of equal size, the reason is simply that 83 is a prime number; a proof would enumerate all of the integers smaller than 83 (the potential group sizes) and show that no number produced a full set of groups without remainder; such a proof would serve as a reason that explains to the students (for example) why some groups were smaller than others.

Not all reasons are equal, however. As Gowers points out in the case of the normality of $\pi$, it is not just that we do not expect there to be any reason for $\pi$ to be not normal: we do not expect there to be a ``good'' reason. To see this in more detail, let us assume, for the moment, that Goldbach's Conjecture \emph{is}, in fact, false, and that there \emph{is} some even number, $n$, a thousand digits long (say), which violates Goldbach's Conjecture. Were this to be the case, we would have a counter-GC theorem that said $\neg \exists x,y~\mathrm{s.t.}~x, y \in\mathrm{Primes}~\&~(x + y) = n$. 

One way to prove this theorem---guaranteed to work---would be to list all primes $p$ less than $n$, and show that in each case $(n-p)$ is not prime. Such a proof would be easy to produce---any undergraduate could write the code to do so, even in a formal system such as Lean---and would, in some sense, provide a reason for why the GC is false. 

Such a reason, however, would be transparently unsatisyfing; whatever other faults it has, the enumeration proof of the counter-GC is extraordinarily long. It is, in the language of computational complexity, \emph{exponentially} (or at least superpolynomially) long in the length of the original statement: the number $n$ has a decimal expansion in (approximately) $\log_2{n}$ bits, but the proof itself has at least $n/\log{n}$ clauses, one for each of the primes less than $n$. This is where Gower's no-coincidence principle gets its power: not that there would exist no reason for counter-GC to be true, but that other considerations suggest that there would be no \emph{good} reason.\footnote{The existence of a very large counterexample would not necessarily imply that there were no good reasons. Imagine that we were told (by an advanced civilization, say) that the GC failed for a number $n$ equal to the Ramsey number $R(6,6)$. In this case, we would rejoice---we would feel that there \emph{has} to be some good reason for that coincidence, surely, and a corresponding proof that was not $R(6,6)/\log{R(6,6)}$ clauses long.}

Listing the sufficient conditions for a proof to be a good reason is difficult; not all succinct proofs are proofs ``from the book''~\cite{aigner2018proofs}, and a proof can fail to be satisfying, illuminating, or otherwise good in the sense Gowers intends for all sorts of reasons. But the badness of an exponentially-long proof is hard to ignore. We can propose a weaker version of Gowers' original principle, the {\bf weak no-coincidence principle}: ``if an apparently outrageous coincidence happens in mathematics, then there is a reasonably efficient proof for it''.

\section{Proof Complexity and the Unreasonable Effectiveness of Mathematics in Mathematics}
\label{pc}

The authority for Gowers' claim, and our weaker restatement of it, rests upon the testimonies of practicing mathematicians. Taking it at face value, we can ask how ``unreasonable'' is this state of affairs, where we can reason backwards from the hypothesised absence of a good reason to the truth of the original statement? To answer this question requires a brief detour into an area of computational complexity known as proof complexity, and the theory of propositional proof systems. 

We can think of a proof, $X$, of a (true) proposition $A$ as something that helps one see that $A$ follows from the axioms of the theory. A propositional proof system (PPS) is an idealized model of this process: it is an algorithm $P$ that outputs tautologies---\emph{all} tautologies (\emph{i.e.}, the proof system is complete) and \emph{only} tautologies (\emph{i.e.}, the proof system is sound). For every tautology $A$, there is a proof $X$ such that $P(X)$ outputs $A$, and (crucially) the algorithm runs in polynomial time in the length of $X$. That $P$ runs in polynomial time captures the idea that a proof is an aid to discovering what is true: rather than (for example) forcing the reader to conduct an exponentially-long brute-force search over all ways to combine the axioms, $X$ provides a road-map for cutting an efficient path to the truth you care about.

A PPS may be an idealization of mathematical practice, but it is a useful one. Running a system such as Lean on a personal computer, for example, approximates a propositional proof system. Errors of judgement aside, it also provides a reasonable approximation of the proof-making, and proof-reading, activity of ordinary human mathematicians. Consider, for example, the standard lacunae, sometimes marked by phrases such as ``obviously'', ``as can be easily seen'', and so forth. Under the assumption that these can be filled in without a significant increase (exponential blow-up) in the search-time required by the reader, such proofs are valid parts of the informal proof systems that constitute written mathematics.

Proof systems promise that a proof, should it be found, can be efficiently checked, but they say nothing about the length of the proof itself. A logical system may well contain many sequences of ``hard'' tautologies, where the length of the proof $X$ for each tautology grows superpolynomially with the length of the tautology $A$. Indeed, a key result in proof complexity,\footnote{As with many results in computational complexity, these rely upon widely believed---but unproven---assumptions about the complexity hierarchy; in this case, that $\NP\neq\coNP$.} (due to Sch{\" o}ning~\cite{hard_proofs}; see Ref.~\cite{monroe2022average} for recent progress) shows that such sequences are dense: if we consider all tautologies of length $n$ (in some binary encoding), then the number of tautologies that are hard in this fashion has density at least $2^{\epsilon n}$. Follow-up results~\cite{instance_hard} show that these hard tautologies can not be eliminated by clever use of axioms---short of adding an axiom of similar length for each  tautology to be proved, hard-to-prove tautologies remain.

The results provide a formal version for Gowers' suggestion, at the end of his paper, that mathematics suffers (or, rather, benefits) from a problem of unreasonable effectiveness. The no-coincidence principle says that true statements have good reasons, but there are a large number of mathematical claims that are (1) tautologies, \emph{i.e.}, true, but that (2) have no efficient proofs. If we accept, for the time being, a correspondence between good reasons and efficient proofs, it is hard to see how mathematicians, in the process of constructing proofs, should consistently contemplate truths that happen to have efficient proofs. Hard-to-prove tautologies not only abound, but are hard to avoid: meta-complexity results suggest that there is no efficient algorithm for showing that a tautology is ``one of the bad ones'', \emph{i.e.}, one that violates the no-coincidence principle~\cite{meta}. 

\section{Do good reasons imply good proofs?}

One way out of this impasse is to break the correspondence between good reasons and efficient proofs: to say, in other words, that there are mathematical truths that have good reasons to be true, but do not have efficient proofs.\footnote{There is, of course, a ``zeroth'' way out: a strong anthropocentrism, in the style of Penrose~\cite{penrose1999emperor}, where there is something strange about human cognition that enables it to circumvent the basic no-go results of the previous section.}

It is certainly conceivable that there are good reasons for things that are not efficiently provable. The Riemann Hypothesis (RH) itself might provide just such an example. As Gowers points out, RH is known to imply a large number of unexpected consequences in other, apparently unrelated, fields~\cite{sarnak}, and this might be seen to provide a reason for it to be true---even though such reasons, being abductive rather than deductive, provide no promise of an efficient proof from first principles. This pattern is common in the empirical sciences, and especially in theoretical physics, where it goes by the name consilience~\cite{WOJTOWICZ2020981}, and where ``unificatory'' theories---\emph{i.e.}, ones that, in describing a fundamental force, can be shown to imply the existence of other observed, but conceptually unrelated, forces as a consequence---are highly prized.

Breaking the reason-proof correspondence would have the unfortunate effect of making the no-coincidence principle a dangerous guide: mathematicians use their intuitions about what is true to guide what to prove next, but---with the exception of grand claims such as the Riemann Hypothesis or $\P\neq\NP$ where conditional results are prized---it is of little use for a mathematician to believe a lemma she can not (efficiently) prove. In this sense, the empirical success of the no-coincidence principle as a guide for getting along with the job of proving things provides evidence in favor of the correspondence.

\section{What if the No-Coincidence Principle were false?}
\label{trials}

A second way out of this impasse is to consider the possibility that the no-coincidence principle may be true often enough to serve as a heuristic for the practice of mathematics, but false in general. If this were the case, it would provide only an imperfect guide to which mathematical statements were true, and a mathematician who relied on it would (1) occasionally try to prove false things, and (2) occasionally avoid relying on true things.

In cases of the first type, a mathematician might encounter a claim $C$, whose falsehood would imply a coincidence without a good reason. Claim $C$ would, in fact, be false, but the mathematician would have no way of discovering that fact, because (by assumption) no efficient proof for it would exist. Assuming some sort of time-limit on her reasoning, she would, eventually, give up, with perhaps a feeling that, had she been more able, a proof might have been found.

In the second type, a mathematician would encounter a (true) claim $C^\prime$, whose truth would imply an extreme coincidence without a good reason. Under these circumstances, she would discard it in her proof-seeking strategy. Discarding such a truth, however, would, counterintuitively, be a benefit to her, since the absence of a good reason would (assuming the correspondence between good reasons and efficient proofs) imply the absence of an efficient proof, and she will have simply saved herself the useless labor.

The testimonial evidence provided by practicing mathematicians, of course, seems to argue against the failure of the no-coincidence principle: when we encounter claims of the first type, we very often find that they're true.

That said, there are (apparent) counter-examples. The prime-counting function, $\pi(n)$, which counts the number of primes less than or equal to $n$, is smaller than the logarithmic integral function, $\mathrm{li}(n)$---at least for small $n$. However, J.E. Littlewood proved (in 1914) that $\mathrm{li}(n)$ must, at some point, become larger than $\pi(n)$, and that the two functions in fact cross each other an infinite number of times. The first point at which the functions cross is likely to be extraordinarily large; our current best upper bound is (assuming the Riemann Hypothesis) $e^{e^{e^{79}}}$, a number with more than $10^{10^{33}}$ digits. Under a Gowers-like quasi-Bayesian analysis, the particular value of this crossing point is necessarily extraordinarily unlikely (a great coincidence) and thus ought to have a very good and (in Gower's analysis) easily-found reason for taking the value it does---but both the number, and its reason for being so, has eluded us for more than a century.

Whether or not the crossing-point example is a true counter-example, it is also notable that much testimonial evidence, perhaps necessarily, focuses on famous examples and succinctly phrased claims. These situations are far removed from the actual practice of doing mathematics. It may well be the case that mathematics ``in the weeds'', as it unfolds minute-by-minute in someone's office involves multiple encounters with things that ought to be true but (from the mathematician's perspective) seem to be unprovable with the resources to hand. 

This raises the possibility that a self-fulfilling prophecy is in operation: use of the no-coincidence principle may guide us through the thickets, at the cost of hiding truths.\footnote{Such an outcome would need, at the very least, to be an imperfect guide, or one that depended sensitively on the particular forms that propositions took, lest it fall afoul of the meta-complexity results mentioned at the end of Section~\ref{pc}.} There may well be vast areas of mathematics that we are guided away from, because we refuse to entertain the truth of unreasonable lemmas.

\section{Negotiating the Field of Efficient Proofs}

Let us, for the moment, accept the results of proof complexity: there are many mathematical statements that do not have efficient proofs but that are, nonetheless, true. The no-coincidence principle is then simply a useful heuristic for getting mathematics done, a useful untruth that dissuades mathematicians from investing too much effort in statements that might be true, but whose proofs would be (truly) unwieldy. 

The demotion of the no-coincidence principle to the no-coincidence heuristic may be help explain something sometimes remarked-upon in the philosophy and history of mathematics: a distinct shift in the nature of mathematical reasoning from\footnote{In Ref.~\cite{zala}'s taxonomy; see also the earlier Ref.~\cite{laut}'s distinction between ``classical'' and ``contemporary''.} the ``classical'' era (mid-17th to mid-19th Century) to the ``modern'' (mid-19th to mid-20th Century) and ``contemporary'' eras (mid-20th Century to present day).

In this view, the evolution of mathematics is characterised by a new sense of what, exactly, mathematics is doing. In the classical era, mathematics might profitably be viewed as the task of proving facts about mathematical objects whose rough outlines are reasonably stable and largely given ahead of time. Such an ``expressivist''~\cite{ab_ex} view is particularly natural in examples where mathematics had close ties to the physical sciences: in the 18th Century, when Lagrange discovered the divergence theorem, there was a natural separation between, on the one hand, his proof and, on the other, the object (continuous vector fields in real-valued spaces) that the proof is establishing facts about.

In the contemporary era, by contrast, proofs and the objects they are about tend to co-vary to a much greater extent. Ref.~\cite{zala} provides testimony from the mathematician Alexander Grothendieck to this effect:
\begin{quote}
    There is a continual and uninterrupted back-and-forth movement here, between the \emph{apprehension} of things and the \emph{expression} of what has been apprehended, by way of a language that has been refined and recreated as the work unwinds, under the constant pressure of immediate needs. (Alexander Grothendieck, \emph{Récoltes et Semailles}, 1985, quoted in~\cite{zala}, 153).
\end{quote}
In the language of this paper, we can make sense of Grothendieck's testimony by seeing mathematicians as negotiating the field of what is provable: in the process of discovery, objects of inquiry are modified until they are such that interesting facts about them can be efficiently proved. This instability plays out against the pressure of immediate needs, among which there is the mathematician's need to prove, as well as intuit, what is in her sights. In contrast to the classical mathematician presented with a question about an object, contemporary mathematics is far more concerned with the discovery of questions to ask, and the invention of proto-objects to ask them about. 

Such a view is consistent with Ref.~\cite{zala}'s account of the ``transitory ontology'' of contemporary mathematics---an account built on both testimony from contemporary mathematicians such as Grothendieck, Lawvere, and Lax, and a reading of the history of research programs in contemporary mathematics, where mathematical objects ``cease to be fixed, stable, classical and well founded ... and tend toward the mobile, the unstable, the nonclassical, and the merely contextually founded''.

Seen in this way, the apparent anthropocentrism that enables mathematical progress takes a different light: the human discipline of (contemporary) mathematics is the practice of making a human-cognizable home within a larger, and mostly hostile, territory. In this proof-first view, much---though not all---of the unreasonableness dissolves. This dissolution comes at a cost, because the practice of mathematics can no longer be seen to be ``about'' stable, or slowly evolving, objects (sets, categories, types, and so on) at all.

\section{Alien Lemmas}

Exponentially long proofs are exponentially long no matter who tries to construct them, but it is natural to ask what would happen if mathematicians were to turn to the far more patient tools provided by automated proof systems such as Coq or Lean. What strange new theorems we might prove, the sanguine might ask, if we only burned enough computer power to establish the unpalatable, alien lemmas they depend upon!

The possibility of relying on mathematical claims without the insight provided by an efficient proof is a strange one,\footnote{See, \emph{e.g.}, Ludwig Wittgenstein's remark: ``We might say, `A mathematical proposition is a pointer to an insight'. The assumption that no insight corresponded to it would reduce it to utter nonsense.''~\cite{phil}, 212.} but all heuristics must fail in some limit, and the failure of the no-coincidence theorem opens the possibility of both alien lemmas, and new forms of mathematics that might be founded on them. Future mathematicians may, indeed, rely on the power of automated assistants akin to the ``AlephZero'' conjectured by Akshay Venkatesh~\cite{aleph}, to discover all kinds of uncanny, unreasonable, ``glitches''~\cite{dedeo}. The process might appear similar to how mathematicians can achieve new results in disparate fields by assuming the Riemann Hypothesis---with a key distinction that these alien lemmas would be known, rather than simply assumed, to be true.

There are reasons to temper one's enthusiasm for such a project. The unwieldy proofs that automated search might uncover are not only exponentially long. They also---as shown by Ref.~\cite{instance_hard}---lack the internal structure that would enable efficient compression, and are thus resistant to being simplified, in a Grothendieck-like fashion, through the creative elaboration of novel concepts. This lack of structure may also make them poor candidates for the analogy-making that often drives mathematical progress~\cite{cor}. Alien lemmas might open up new fields, but their proofs would provide no clue for how to explore them.

Searching for such lemmas is not an easy task, if the search is in aid of producing generative truths. It is easy (given background assumptions) to write down hard-to-prove tautologies that take the form of apparently arbitrary claims about this or that hard-to-invert one-way function, with no clear downstream utility. A similar phenomenon occurs in the case of undecidability in ZFC: while certain questions of extreme interest (\emph{e.g.}, the continuum hypothesis) are undecidable, the DPRM Theorem says that there are plenty of \emph{prima facie} boring ones, concerning long and apparently structureless diophantine equations, as well. A converse question is how common it is, in actual mathematical practice, for a mathematician to stumble upon (and then unknowingly deny, per Section~\ref{trials}) hard-to-prove truths. A recent near-miss in this domain might well be the recent result concerning the bunkbed conjecture~\cite{gladkov2024bunkbed}---a near-miss because the authors were, in fact, able to find an efficient formal argument.

\section{Conclusion}

Mathematics is a paradigmatically deductive project: by whatever means its truths are first intuited, they must, in the final analysis, be established by deductive proof. It is also, as mathematicians themselves testify, a paradigmatically reasonable project: mathematical truths are often true for good reasons and this guides the process of proof-making itself. This work has presented a minimal account of what that reasonableness might be. It has then drawn on basic results in computational complexity to show what the combination of deduction and reasonableness implies about the larger architecture of that project.

\vspace{0.5cm} 

{\bf Acknowledgements}. I thank the participants of the First ILIAD Conference in Theoretical Alignment, Berkeley, CA (2024), Paul Livingston, Colin Allen, and Gülce Kardeş for helpful discussions. This work was supported in part by the Survival and Flourishing Fund.

\bibliographystyle{unsrt}  
\bibliography{templateArxiv}

\end{document}